\documentclass[a4paper,11pt]{article}

\usepackage{amsmath,amsfonts,amssymb,amsthm,amscd}
\usepackage{graphicx,mathtools}
\usepackage{perpage}
\usepackage{url}
\usepackage{color}
\usepackage[T1]{fontenc}
\usepackage{hyperref}
\usepackage[a4paper,scale={0.72,0.74},marginratio={1:1},footskip=7mm,headsep=10mm]{geometry}
\usepackage[latin9]{inputenc}
\usepackage{mathrsfs}
\usepackage{tikz}
\usepackage{caption}
\usepackage{subcaption}
\usepackage{float}
\usepackage{enumerate}
\usepackage[english]{babel}
\usepackage{hyperref}
\usepackage{setspace}
\usepackage{rotating}
\usepackage{bbm}
\usepackage{booktabs}
\usepackage{dsfont}

\makeatletter
\def\@captype{figure}
\makeatother

\numberwithin{equation}{section}

\DeclareRobustCommand{\DIEP}{\ensuremath{
{\includegraphics[height=2ex]{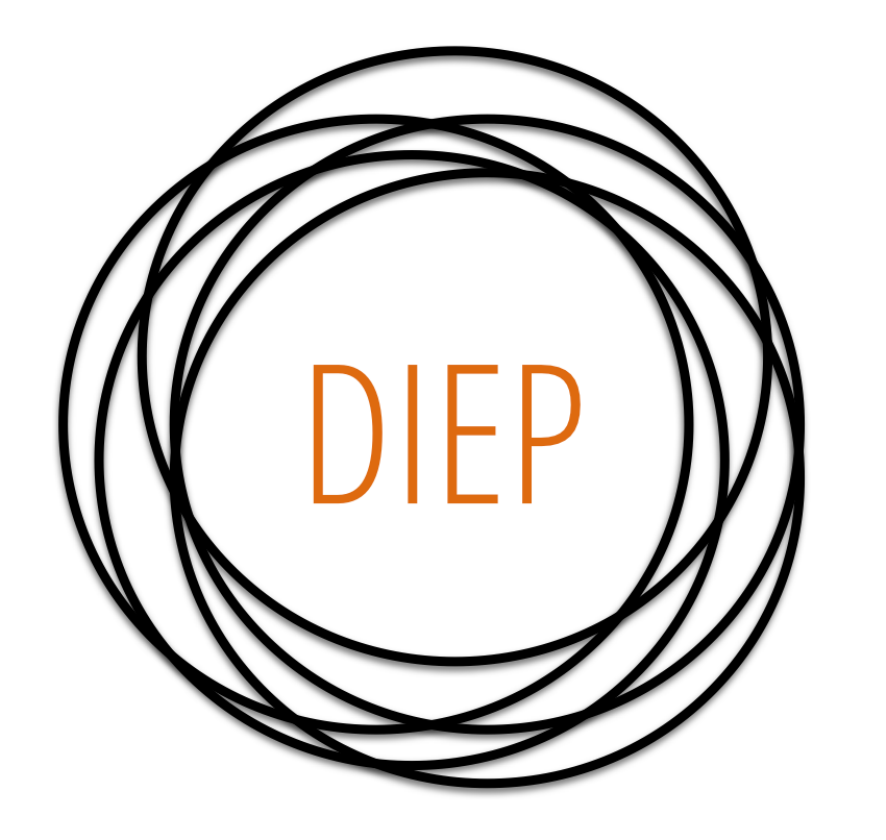}}
}}

\title{Limiting dynamics for Q-learning with memory one in two-player, two-action games}
\author{

J.\ M.\ Meylahn
\footnotemark[12]\thanks{Contact: j.m.meylahn[at]utwente.nl}
\\

L.\ Janssen
\footnotemark[34]
}

\footnotetext[1]{
Department of Applied Mathematics, University of
Twente, The Netherlands.
}

\footnotetext[2]{
\DIEP Dutch Institute of Emergent Phenomena, University of
Amsterdam, The Netherlands.
}

\footnotetext[3]{
Korteweg-de Vries Institute for Mathematics, University of
Amsterdam, The Netherlands.
}

\footnotetext[4]{
Institute for Informatics, University of
Amsterdam, The Netherlands.
}

\begin{document}
\maketitle
\begin{abstract}
We develop a computational method to identify all pure strategy equilibrium points in the strategy space of two-player, two-action repeated games played by Q-learners with one-period memory. This is done by constructing a graph of pure strategy mutual best-responses. We apply this method to the iterated prisoner's dilemma, stag hunt and hawk-dove games. For each of the games we are able to identify all possible equilibrium strategy pairs and the conditions for their existence. We confirm that the mutual best-response dynamics are realized by a Sample Batch Q-learning algorithm in the infinite batch size limit. 

\medskip\noindent
{\it Acknowledgment. }  This work is supported by the Dutch Institute for Emergent Phenomena (DIEP) cluster at the University of Amsterdam.
\end{abstract}

\newpage
\section{Introduction}

As algorithms are increasingly employed in a variety of settings, more situations will arise in which algorithms interact. In the case of (reinforcement) learning algorithms, what this interaction will lead to is not easily predictable. This is due to the algorithms being designed for single-agent settings (or stationary environments). The theoretical foundations for reinforcement learning algorithms in multi-agent settings are thin, with few convergence results \cite{Busoniu2008}.

A particularly relevant application domain in which understanding the dynamics of multi-agent learning is important, is that of algorithmic pricing \cite{Meylahn2022}. Here multiple algorithm learn optimal prices to maximize profits for their respective firms in a market. If these algorithms converge to higher than competitive prices (and thus learning to cooperate or collude), this has detrimental effects for consumers. Recent work has shown that Q-learning with memory can learn to collude in this setting \cite{Calvano2020, Calvano2020b, Klein2021}. 

This article aims at contributing to the theoretical foundations in the case of two-player, two-action games in two ways: (1) by developing a computational method for determining limiting dynamics of Q-learning algorithms with one-period memory and (2) by constructing a representation of deterministic one-period memory Q-learning dynamics. We demonstrate the methods using the iterated prisoner's dilemma (PD) \cite{Rapoport1965}, the stag hunt (SH) game \cite{Rousseau1754} and the hawk-dove (HD) game \cite{Smith1973}.    

The PD is one of the paradigmatic model for studying the emergence of cooperation in social dilemma's \cite{Macy2002}, pricing duopolies \cite{Waltman2008} and other economics games \cite{Rabin1993}. The model's simplicity allows for rigorous analysis, giving insights that can be heuristically extended to more realistic settings. Two simple extensions to the basic model which make significant steps toward becoming realistic while keeping the model tractable are to include a basic learning mechanism for the players and to give the players a period one memory.  

In realistic games, the players of the game have to learn their optimal strategy by learning what the payoffs of the game are and how their opponent is behaving.  Games in which both players are learning fall into the field of multi-agent learning (MAL), in which much recent work has focused on dealing with the learning pathologies that typically arise in MAL (see\cite{Palmer2020} and references therein). For surveys of MAL see  \cite{Busoniu2010, Hernandez2019, Panait2005, Zhang2021}  

One of the distinctions that can be made in the MAL literature is that between single-state and multi-state settings. A multi-state setting can be thought of as a situation where multiple games are being played, with probabilistic transitions occurring between the games \cite{Hennes2008}. One of the most simple multi-state settings is the one where the players have a period one memory. The state the game is in is just the action pair played in the previous round, which both players can observe. This is economically relevant since many of the classical strategies developed for the PD require a one-period memory \cite{Axelrod1981}. The Tit-For-Tat strategy, for example, imitates the action of the opponent in the previous round. 

Among the most studied reinforcement learning algorithms for MAL is Q-learning \cite{Watkins1992}. This is an algorithm designed for Markov decision processes for single agents. In multi-agent settings, convergence of Q-learning is typically poor due to the many pathologies that can occur \cite{Palmer2020}. A promising method for overcoming these pathologies by dealing with the non-stationarity of MAL is to learn in batches. Both agents keep their strategy fixed during the batch in order to collect data in a stationary environment and then simultaneously update their Q-values at the end of the batch. In \cite{Barfuss2019} the authors show that this leads to deterministic dynamics in the infinite batch size limit. The limiting dynamics are a good approximation for large, but finite batch sizes.  

In this paper, we develop a computational method for determining the absorbing states of batch Q-learning with period one memory in two-player, two-action iterated games. These games have been classified for Q-learning without memory in \cite{Wunder2010}. The method allows us to plot what we call a iterated mutual best-response networks (IBRN), which gives the strategy space dynamics in the infinite batch size limit. We show that the absorbing states identified in \cite{Usui2021} for the PD are the \emph{only} absorbing states. In addition to the absorbing sates, we find that limit cycles are possible in the IBRN graph. This is an artifact of the best-response dynamics and should not persist in algorithms with decreasing learning rate.  

The absorbing states identified in this way for batch Q-learning must necessarily be the same as the absorbing states for \emph{any successful} (in the sense of solving the Bellman equations) Q-learning algorithm that converges to a pure strategy (for example by having a decreasing exploration rate).

In Section \ref{sec:setting} we identify notation, define the model and introduce Q-learning. In Section \ref{sec:method} we describe the methods we develop. In Section \ref{sec:results}  and \ref{sec:numerics} we give the theoretical and numerical results respectively. Finally, in Section \ref{sec:conclusion} we conclude and discuss future research.

\section{Setting}
\label{sec:setting}
\subsection{Model}
We consider the same setting as \cite{Usui2021} and we will employ broadly similar notation for ease of reference. The setting is that of an PD played by two Q-learners, labeled $a\in \{1, 2\}$ and using the convention $-a=\{1, 2\}\backslash a$, with one-period memory. Each Q-learner has a choice between two actions: defect (D) or cooperate (C). We denote the actions by $\sigma_{a}\in \{C, D\}$. Due to the one-period memory the game becomes a multi-state game with the possible states being $\boldsymbol{\sigma}=(\sigma_{1}, \sigma_{2})\in \{(D, D), (D,C), (C,D),(C, C)\}$, where the first component of the vector denotes the action of player one and the second component denotes the action of player two. The payoff (or reward) for player $a$ when the game is in state $\boldsymbol{\sigma}$ is denoted by $r_{a}(\boldsymbol{\sigma})$ and takes the following form:
\begin{equation}
\begin{pmatrix}
r_{1}(1), r_{2}(1)  & r_{1}(2), r_{2}(2)\\
r_{1}(3), r_{2}(3) & r_{1}(4), r_{2}(4)
\end{pmatrix}
=
\begin{pmatrix}
p, p  & t, s\\
s, t & r, r
\end{pmatrix},
\end{equation}
where we have enumerated the states as $(D, D)\rightarrow 1, (D,C)\rightarrow 2, (C,D)\rightarrow 3$, and $(C, C)\rightarrow 4$. In addition, we make the standard assumptions that $t>r>p>s$ for the PD, $r>t>p>s$ for the SH game and $t>r>s>p$ for the HD game.

\subsection{Learning}

Each player chooses their actions according to a strategy $\pi_{a}$, which in the Q-learning with period one memory setting is the conditional probability of playing action $\sigma_{a}$ given that the game is in state $\boldsymbol{\sigma}'$, i.e. $\pi_{a}(\sigma_{a}|\boldsymbol{\sigma}')$. 

The strategies of the players are updated through a reinforcement learning process. We specifically consider Q-learning type algorithms, which are designed to learn the actions that maximize the discounted future rewards, i.e. 
\begin{equation}
\label{eq:discountedprofit}
\sum_{n=1}^{\infty}\delta^{n} r_{a}(\boldsymbol{\sigma}(n)),
\end{equation}   
with $0<\delta<1$ being the discount factor and $n$ denoting the time. 

It has been shown that maximizing \eqref{eq:discountedprofit} can be achieved by solving 
\begin{equation}
\label{eq:qfunction}
Q_{a, \boldsymbol{\sigma}(t), \sigma_{a}(t)} = \mathbb{E}[r_{a}(\boldsymbol{\sigma}(t+1))| \boldsymbol{\sigma}(t), \sigma_{a}(t)] + \delta \mathbb{E}\Big[\max_{i\in\{C, D\}} \{Q_{a, \boldsymbol{\sigma}(t+1), i}\} |\ \boldsymbol{\sigma}(t), \sigma_{a}(t)\Big].
\end{equation}
Note that since the actions taken in the previous round determine the current state, we have used the notation that $\boldsymbol{\sigma}(t) = (\sigma_{1}(t-1), \sigma_{2}(t-1))$. The first term on the right-hand side is the expected reward from taking action $\sigma_{a}(t)$ while in state $\boldsymbol{\sigma}(t)$ and the second term is the discounted expected maximum reward of behaving optimally (according to the current estimate) in the state $\boldsymbol{\sigma}(t+1)$ reached by taking action $\sigma_{a}(t)$ in state $\boldsymbol{\sigma}(t)$. The Q-function defined by \eqref{eq:qfunction} is related to Bellman's value function by $V_{a}(\boldsymbol{\sigma}) = \max_{i\in\{C, D\}}\{Q_{a, \boldsymbol{\sigma}, i}\}$.

The Q-learner $a$ solves \eqref{eq:qfunction} by initializing a Q-matrix
\begin{equation}
\begin{pmatrix}
Q_{a, 1,D} & Q_{a, 1,C} \\
Q_{a, 2,D} & Q_{a, 2,C} \\
Q_{a, 3,D} & Q_{a, 3,C}  \\
Q_{a, 4,D} & Q_{a, 4,C}
\end{pmatrix}
\end{equation}
and updating the entries according to
\begin{equation}
\label{eq:Qupdates}
Q_{a, \boldsymbol{\sigma}, \sigma_{a}}(t+1) = (1-\alpha(t))Q_{a, \boldsymbol{\sigma}, \sigma_{a}}(t) + \alpha(t)\Big[r_{a}(\boldsymbol{\sigma}(t+1)) +\delta \max_{i\in\{C, D\}} \{Q_{a, \boldsymbol{\sigma}(t+1), i}\}\Big],
\end{equation}
where 
\begin{equation}
\alpha(t)=
\begin{cases}
\alpha \quad \text{if} \; (\boldsymbol{\sigma}, \sigma_{a}) = (\boldsymbol{\sigma}(t), \sigma_{a}(t))\\
0 \quad \text{otherwise}
\end{cases}
\end{equation}
is the learning rate. Note that $r_{a}(\boldsymbol{\sigma}(t+1))$ is the reward obtained in period $t$ (i.e. the reward of taking action $\sigma_{a}(t)$ while in state $\boldsymbol{\sigma}(t)$, which leads to state $\boldsymbol{\sigma}(t+1)$). In stationary environments the iterative process defined by \eqref{eq:Qupdates} has been shown to converge under some general conditions, for example, in \cite{Jaakkola1994}.

In contrast to the general multi-state setting, the transitions between states here happen deterministically given the actions of the players. This means that the only stochasticity comes from the action-selection mechanism, i.e. the strategies. Such an action-selection mechanism must take into account the exploration-exploitation trade-off, where exploratory actions are conducive to learning and exploitatory actions maximize the reward given the agent's current state of knowledge. The methods used in Q-learning for this balancing act can be split into two main categories: the $\epsilon$-greedy mechanisms and the softmax or Boltzman mechanisms. We will focus on mechanisms of the first type and leave the latter mechanisms for future work. We further distinguish between $\epsilon$-greedy mechanisms that either have a constant $\epsilon$ or a time dependent $\epsilon(n)$, which decreases in time.

With an $\epsilon$-greedy action-selection mechanism, the player chooses the action which it considers to be maximizing its reward with probability $1-\epsilon(n)$ and chooses an action uniformly at random with probability $\epsilon(n)$.  

Q-learning converges to a solution of the Bellman equations in stationary environments. A Q-learner in a multi-agent setting is however facing a non-stationary environment since the strategy of the opponent changes in time. In the case of a decreasing $\epsilon(n)$, both agents will eventually find themselves in an environment that is stationary, since the opponent's strategy converges to a deterministic (pure) strategy.

In order to identify equilibrium points of multi-agent Q-learning, it thus makes sense to study the dynamics of Q-learning type algorithms that are focused on learning pure strategies, since the equilibrium points should be the same. To this end, \cite{Usui2021} conceive of a Q-learning algorithm in which the players alternate between either learning with the use of Q-learning or keeping their strategy fixed. This ensures that the environment is stationary during the learning process for both players. In this way, the players sequentially learn a best-response to a pure strategy.

As an alternative algorithm that has been suggested for multi-agent settings, we consider a batch learning algorithm as in \cite{Barfuss2019}, in which both players act according to a fixed strategy for the duration of the batch and then update their strategies simultaneously. In the infinite batch size limit, this leads to both players simultaneously taking a step towards the best-response given the strategy of the opponent. The absorbing states of this algorithm are the same as the absorbing states of the algorithm in \cite{Usui2021}.

\section{Method}
\label{sec:method}
\subsection{Self-consistent solutions to the Bellman equations}
\label{sec:selfcons}
If we restrict ourselves to pure strategies, each player can pick one of two actions for each of the four states, leading to a total of $2^{4}$ strategies. As a result, the total number of pure strategy \emph{pairs} is $2^{8}=256$. Using computer algebra software, such as Mathematica, we can automate the calculations in \cite{Usui2021}, performed for the $16$ symmetric strategy pairs, to repeat them for all possible strategies pairs. This allows us to identify a best-response for each strategy.

To this end, we identify each strategy with a 4 dimensional vector. The four entries encode the strategy. A zero in the first entry means that $Q_{1, 1, D}>Q_{1, 1, C}$, while a one means that $Q_{1, 1, D}<Q_{1, 1, C}$. The vector $(0, 0, 0, 0)$, for example, represents the all defect (All-D) strategy. The vector $(1, 0, 0, 1)$ in contrast represents the Win-Stay-Lose-Shift (WSLS) strategy. In a similar fashion, we define a strategy pair to be an 8 dimensional vector in which the first four entries encode the strategy of player 1 and the last four entries encode the strategy of player 2.  

Each strategy pair vector leads to a system of sixteen linear equations that can be solved simultaneously for the sixteen Q-values. The Q-values are then expressed in terms of the model parameters $t, r, p$ and $s$. The vector however also encodes an assumption on the inequalities obtaining between the Q-values. We then perform a self-consistency check to see if the resulting Q-values satisfy the assumed inequalities given the assumptions on the model parameters. In this way, we can determine under which conditions there is a valid solution to the Bellman equations given the strategy pair.

To illustrate the computation, we will consider the PD where both players are playing the all defect strategy, i.e., the vector $(0, 0, 0, 0, 0, 0, 0, 0)$.
The Bellman equations \eqref{eq:qfunction} for player 1 can be written in terms of indicator functions when considering pure strategies. The equations for player 1 become
\begin{align}
Q_{1, 1, D} =& \mathds{1}_{\{Q_{2, 1, D}>Q_{2, 1, C}\}}(p+\delta Q_{1, 1, D}\mathds{1}_{\{Q_{1, 1, D}>Q_{1, 1, C}\}}+\delta Q_{1, 1, C}\mathds{1}_{\{Q_{1, 1, D}<Q_{1, 1, C}\}})\nonumber\\
&+ \mathds{1}_{\{Q_{2, 1, D}<Q_{2, 1, C}\}}(t+\delta Q_{1, 2, D}\mathds{1}_{\{Q_{1, 2, D}>Q_{1, 2, C}\}}+\delta Q_{1, 2, C}\mathds{1}_{\{Q_{1, 2, D}<Q_{1, 2, C}\}})\\
Q_{1, 1, C} =& \mathds{1}_{\{Q_{2, 1, D}>Q_{2, 1, C}\}}(s+\delta Q_{1, 3, D}\mathds{1}_{\{Q_{1, 3, D}>Q_{1, 3, C}\}}+\delta Q_{1, 3, C}\mathds{1}_{\{Q_{1, 3, D}<Q_{1, 3, C}\}})\nonumber\\
&+ \mathds{1}_{\{Q_{2, 1, D}<Q_{2, 1, C}\}}(r+\delta Q_{1, 4, D}\mathds{1}_{\{Q_{1, 4, D}>Q_{1, 4, C}\}}+\delta Q_{1, 4, C}\mathds{1}_{\{Q_{1, 4, D}<Q_{1, 4, C}\}})\\
Q_{1, 2, D} =& \mathds{1}_{\{Q_{2, 2, D}>Q_{2, 2, C}\}}(p+\delta Q_{1, 1, D}\mathds{1}_{\{Q_{1, 1, D}>Q_{1, 1, C}\}}+\delta Q_{1, 1, C}\mathds{1}_{\{Q_{1, 1, D}<Q_{1, 1, C}\}})\nonumber\\
&+ \mathds{1}_{\{Q_{2, 2, D}<Q_{2, 2, C}\}}(t+\delta Q_{1, 2, D}\mathds{1}_{\{Q_{1, 2, D}>Q_{1, 2, C}\}}+\delta Q_{1, 2, C}\mathds{1}_{\{Q_{1, 2, D}<Q_{1, 2, C}\}})\\
Q_{1, 2, C} =& \mathds{1}_{\{Q_{2, 2, D}>Q_{2, 2, C}\}}(s+\delta Q_{1, 3, D}\mathds{1}_{\{Q_{1, 3, D}>Q_{1, 3, C}\}}+\delta Q_{1, 3, C}\mathds{1}_{\{Q_{1, 3, D}<Q_{1, 3, C}\}})\nonumber\\
&+ \mathds{1}_{\{Q_{2, 2, D}<Q_{2, 2, C}\}}(r+\delta Q_{1, 4, D}\mathds{1}_{\{Q_{1, 4, D}>Q_{1, 4, C}\}}+\delta Q_{1, 4, C}\mathds{1}_{\{Q_{1, 4, D}<Q_{1, 4, C}\}})\\
Q_{1, 3, D} =& \mathds{1}_{\{Q_{2, 3, D}>Q_{2, 3, C}\}}(p+\delta Q_{1, 1, D}\mathds{1}_{\{Q_{1, 1, D}>Q_{1, 1, C}\}}+\delta Q_{1, 1, C}\mathds{1}_{\{Q_{1, 1, D}<Q_{1, 1, C}\}})\nonumber\\
&+ \mathds{1}_{\{Q_{2, 3, D}<Q_{2, 3, C}\}}(t+\delta Q_{1, 2, D}\mathds{1}_{\{Q_{1, 2, D}>Q_{1, 2, C}\}}+\delta Q_{1, 2, C}\mathds{1}_{\{Q_{1, 2, D}<Q_{1, 2, C}\}})\\
Q_{1, 3, C} =& \mathds{1}_{\{Q_{2, 3, D}>Q_{2, 3, C}\}}(s+\delta Q_{1, 3, D}\mathds{1}_{\{Q_{1, 3, D}>Q_{1, 3, C}\}}+\delta Q_{1, 3, C}\mathds{1}_{\{Q_{1, 3, D}<Q_{1, 3, C}\}})\nonumber\\
&+ \mathds{1}_{\{Q_{2, 3, D}<Q_{2, 3, C}\}}(r+\delta Q_{1, 4, D}\mathds{1}_{\{Q_{1, 4, D}>Q_{1, 4, C}\}}+\delta Q_{1, 4, C}\mathds{1}_{\{Q_{1, 4, D}<Q_{1, 4, C}\}})\\
Q_{1, 4, D} =& \mathds{1}_{\{Q_{2, 4, D}>Q_{2, 4, C}\}}(p+\delta Q_{1, 1, D}\mathds{1}_{\{Q_{1, 1, D}>Q_{1, 1, C}\}}+\delta Q_{1, 1, C}\mathds{1}_{\{Q_{1, 1, D}<Q_{1, 1, C}\}})\nonumber\\
&+ \mathds{1}_{\{Q_{2, 4, D}<Q_{2, 4, C}\}}(t+\delta Q_{1, 2, D}\mathds{1}_{\{Q_{1, 2, D}>Q_{1, 2, C}\}}+\delta Q_{1, 2, C}\mathds{1}_{\{Q_{1, 2, D}<Q_{1, 2, C}\}})\\
Q_{1, 4, C} =& \mathds{1}_{\{Q_{2, 4, D}>Q_{2, 4, C}\}}(s+\delta Q_{1, 3, D}\mathds{1}_{\{Q_{1, 3, D}>Q_{1, 3, C}\}}+\delta Q_{1, 3, C}\mathds{1}_{\{Q_{1, 3, D}<Q_{1, 3, C}\}})\nonumber\\
&+ \mathds{1}_{\{Q_{2, 4, D}<Q_{2, 4, C}\}}(r+\delta Q_{1, 4, D}\mathds{1}_{\{Q_{1, 4, D}>Q_{1, 4, C}\}}+\delta Q_{1, 4, C}\mathds{1}_{\{Q_{1, 4, D}<Q_{1, 4, C}\}}).
\end{align}
These are reduced to
\begin{align}
\label{eq:BellallD}
Q_{1, 1, D} = p +\delta Q_{1, 1, D}, \; Q_{1, 1, C} = s +\delta Q_{1, 3, D}, \; Q_{1, 2, D} = p +\delta Q_{1, 1, D}, \; Q_{1, 2, C} = s +\delta Q_{1, 3, D},\nonumber\\
Q_{1, 3, D} = p +\delta Q_{1, 1, D}, \; Q_{1, 3, C} = s +\delta Q_{1, 3, D}, \; Q_{1, 4, D} = p +\delta Q_{1, 1, D}, \; Q_{1, 4, C} = s +\delta Q_{1, 3, D}
\end{align}
when considering the strategy pair of both players using the All-D strategy. Solving \eqref{eq:BellallD} leads to 
\begin{align}
&Q_{1, 1, D}=Q_{1, 2, D}=Q_{1, 3, D}=Q_{1, 4, D} = \frac{p}{1-\delta}\\
&Q_{1, 1, C}=Q_{1, 2, C}=Q_{1, 3, C}=Q_{1, 4, C} = s + \frac{\delta p}{1-\delta}.
\end{align}
Since this is a symmetric strategy pair, the solutions to the Bellman equations for player 2 are the same.

To see if this is a valid solution or not, we check that the inequalities of the assumed strategy pair are satisfiable by the model parameters. In this case, the inequalities require that
\begin{equation}
Q_{a, i, D} > Q_{a, i, C} \quad \text{for all } a\in\{1, 2\} \text{ and } i\in \{1, 2, 3, 4\},
\end{equation}    
so that we must have 
\begin{equation}
\frac{p}{1-\delta} > s + \frac{\delta p}{1-\delta},
\end{equation}
or equivalently
\begin{equation}
p(1-\delta) > s(1-\delta). 
\end{equation}
This is always satisfied since $p>s$ for the PD. We conclude that the All-D strategy is the best-response against an opponent playing the All-D strategy.

\subsection{Best-response networks} 
\label{sec:BRN} 
We can now construct two types of direct networks for each choice of the parameters $t, r, p, s$ and $\delta$ from the pure strategy best-responses. The first type we call best-response networks (BRN) and they are constructed as follows. We convert each of the 16 strategies into a number: $(0, 0, 0, 0)\rightarrow 0$, $(0, 0, 0, 1)\rightarrow 1$, etc.\footnote{To do this consistently we use the binary representation of numbers.} These become the labels for the nodes of our graph and we draw a directed edge from each strategy to its best-response. 

The second type of networks we call iterated mutual best-response networks (IBRN). For these we convert the strategy pairs into numbers so that $(0, 0, 0, 0, 0, 0, 0, 0) \rightarrow 0$, $(0, 0, 0, 0, 0, 0, 0, 1) \rightarrow 1$ etc. Now we draw a directed edge from a strategy pair to the strategy pair that contains the best-response to player ones strategy in the last four entries and the best-response to player two's strategy in the first four entries. The IBRN encode the mutual best-response dynamics where in each round each player switches their strategy to the strategy that is a best-response to their opponents previous strategy (see \cite{Fudenberg1991}). 

Absorbing states for the mutual best-response dynamics will appear as nodes with self-loops in the IBRN. These absorbing states correspond to Nash equilibria for the dynamics. Symmetric absorbing states will appear as self-loops in the BRN as well, but non-symmetric absorbing states will appear as reciprocal edges. In addition to the absorbing states, the IBRN will exhibit limit cycles. These will consist of strategy pairs in which each player plays a strategy that is part of an absorbing state, but not the same absorbing state. In this situation, the mutual best-response dynamics dictate that the players will switch to the opponent's strategy at each round. This is an example of the mis-coordination learning pathology discussed in \cite{Palmer2020}. Our focus in this article will be on absorbing states. 

\subsection{Classification of strategy pairs}
\label{sec:classification}
In this section we introduce two classifications that allow us to identify the type of the absorbing strategy pairs. The first classification is based on the type of behavior the strategy pair gives rise to and the second classification is based on the symmetries or anti-symmetries in the strategy pair. 

\subsubsection*{From strategies to actions}
\label{sec:meth:act}
Knowing which strategy pair is being played is still not necessarily informative as to which actions will be taken by the algorithms. As an example, consider the symmetric strategy pair of both players using the Tit-For-Tat strategy with exploration. If the system starts in the CC or DD state, it will remain there until one of the players explores and picks D or C, respectively. Then the system will oscillate between the DC and CD states until one of the players again explores to synchronize their action with that of the opponent. This means that the state graph of the game has three disconnected components. Transitions between components occur only due to exploration.

\begin{figure}
\centering
\includegraphics[scale=0.5]{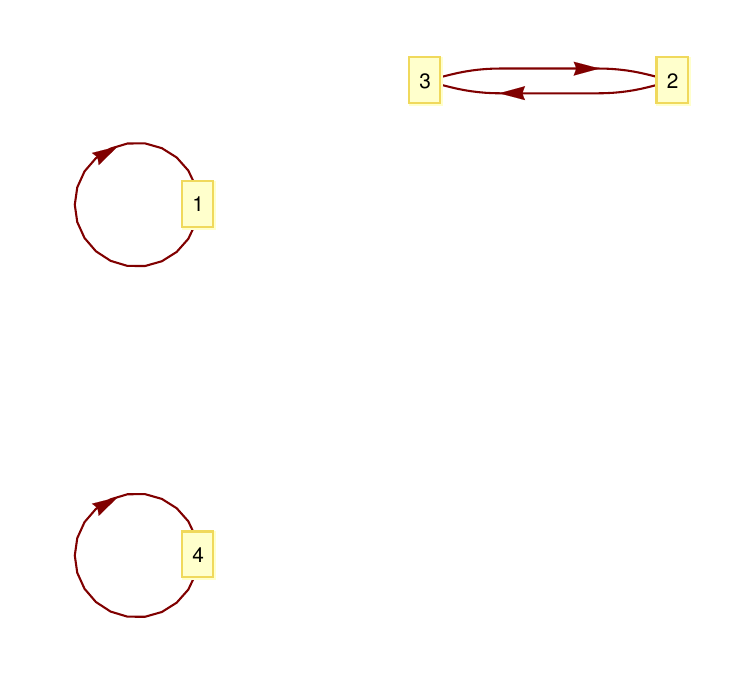}
\caption{Game state transitions of the symmetric TFT pure strategy pair}
\label{fig:TFTpure}
\end{figure}

The vector for this state is given by $(0, 1, 0, 1, 0, 0, 1, 1)$ (note that states 2 and 3 require opposite responses from players 1 and 2). The transitions between states of the game when both players play a pure strategy are shown in Figure \ref{fig:TFTpure}.

We can include the transitions due to exploration as follows. Let the probability that no one explores be $(1-\epsilon)(1-\epsilon) = A$, the probability that player 1 explores be $\epsilon(1-\epsilon)=B$ (this is the same as the probability that player 2 explores) and the probability that both explore be $\epsilon^{2}=C$. Then the transition (stochastic) matrix for the game with the symmetric TFT strategy pair and $\epsilon$-greedy exploration is given by 
\begin{equation}
P = \begin{pmatrix}
A + B + \frac{C}{4} & \frac{B}{2} + \frac{C}{4} &\frac{B}{2} + \frac{C}{4} & \frac{C}{4}\\
\frac{B}{2} + \frac{C}{4} & \frac{C}{4} &  A + B + \frac{C}{4} & \frac{B}{2} + \frac{C}{4} \\
\frac{B}{2} + \frac{C}{4} & A + B + \frac{C}{4} & \frac{C}{4} & \frac{B}{2} + \frac{C}{4}\\
\frac{C}{4} & \frac{B}{2} + \frac{C}{4} &\frac{B}{2} + \frac{C}{4} & A + B + \frac{C}{4}
\end{pmatrix}.
\end{equation}
The eigenvector of $P$ with eigenvalue 1 will give the stationary distribution on the game states and therefore tell us how much time is spent in each state. For this matrix the stationary distribution is the uniform distribution on the states, i.e., $(1/4, 1/4, 1/4, 1/4)$. This means that the game spends equal amounts of time at each state. 

A similar calculation yields the following stationary distribution for the symmetric Grim Trigger strategy pair:
\begin{equation}
\mathcal{N}\begin{pmatrix}
\frac{4-5\epsilon - 2\epsilon^{2}}{\epsilon} & 2-\epsilon & 2 - \epsilon & 1
\end{pmatrix},
\end{equation}
where $\mathcal{N}$ is a normalizing constant. For $\epsilon=0.1$ this is $(0.997, 0.054, 0.054, 0.028)$, which shows that even though coordination is possible in the GT symmetric strategy pair it occurs very little (as evidence by the $0.028$ probability of being in state 4).

By calculating the stationary distributions of the pure strategy pairs, we can identify which states of the game the learners will spend the most time in, giving us an idea of the actions used and the expected payoff. Using this, we can classify the strategy pairs by the resulting stationary distributions. In our case, we identify a strategy pair as being conducive to coordination (CC) or not (NCC). A strategy pair is conducive to coordination when the strategy pair assigns positive probability to being in the $(C, C)$ state. This provides an avenue for studying the likelihood of coordination under pure mutual best-response dynamics.

\subsubsection*{Structural symmetry}
We find that all absorbing states in the three games we consider fall in one of the following four categories.
\begin{description}
    \item[Symmetric (S)] The players play the same action in all states.\\ An example is $(0,1,1,0,0,1,1,0)=102$.
    \item[Complemented middle (CM)] The players play the same action in the states $(0,0)$ and $(1,1)$, but complimentary actions in the states $(0,1)$ and $(1,0)$.\\ An example is $(0,0,1,0,0,1,0,0)=36$.
    \item[Complemented sides (CS)] The players play the same action in the states $(0,1)$ and $(1,0)$, but complimentary actions in the states $(0,0)$ and $(1,1)$.\\ An example is $(0,0,1,1,1,0,1,0)=58$.
    \item[Complemented (C)] The players play opposite actions in all states.\\ An example is $(0,1,1,1,1,0,0,0)=120$.
\end{description}

\section{Results}
\label{sec:results}
In this section we collect the critical conditions at which the BRN change, give examples of the resulting networks and identify the absorbing state of the mutual best-response dynamics. We do this for the prisoner's dilemma, the stag hunt and the hawk-dove games, each of which corresponds to a different ordering of the model parameters $t, r, p$ and $s$. Without loss of generality and for ease of exposition we will always normalize the games by setting $\max\{t, r, p, s\} = 1$ and $\min\{t, r, p, s\} = 0$. This reduces the dimension of the parameter space to 3 and thus allows us to visualize the phase space.

\subsection*{Prisoner's dilemma}
\label{sec:resultsPD}
Using the procedure outlined in Section \ref{sec:selfcons}, we find that the symmetric absorbing states found by \cite{Usui2021} are the only possible absorbing states. This means that there are no non-symmetric strategy pairs that solve the Bellman equations in this setting. 

The only possible pure strategy solutions pairs are the symmetric all defect (All-D) strategy represented by node 0, the Grim Trigger (GT) strategy represented by node 17 and the Win-Stay-Lose-Shift (WSLS) strategy represented by node 153 (the last strategy is also known as the one-period-punishment strategy or the Pavlov strategy). The All-D strategy is always possible, while there are restrictions on the values of the model parameters determining when the other two are possible. This is summarized in Table \ref{tab:absorbingPD}.

\begin{table}[!ht]
    \centering
    \begin{tabular}{c|c|c|c}\toprule
        \textbf{Policy-pair} & \textbf{Type} & \textbf{Conditions} & \textbf{Region}\\
        \hline
        $0$ & S, NCC & Always & $(0)$\\
        $17$ & S, NCC &$r+\delta  > 1 + \delta p$ & $(17)$\\
        $153$ & S, CC & $r + \delta r > 1 + \delta p$ & $(153)$\\
        \bottomrule
    \end{tabular}
    \caption{All strategy pairs solving the Bellman equations self-consistently for the prisoner's dilemma, their type and the condition under which they exist.}
    \label{tab:absorbingPD}
\end{table}

\begin{figure}
\centering
\includegraphics[scale=1.8]{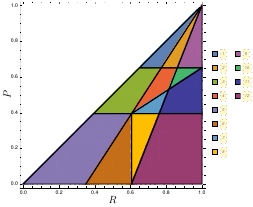}
\includegraphics[scale=1.65]{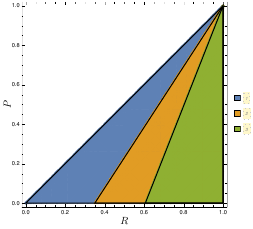}
\caption{Left: Phase diagram for the possible BRNs in the prisoner's dilemma Right: Phase diagram for the possible equilibria in the prisoner's dilemma (only node 0 in region 1, nodes 0 and 1 in region 2 and nodes 0, 1 and 9 in region 3). For both plots we have $s=0, t=1$ and $\delta=0.65$. }
\label{fig:regionplotPDd65}
\end{figure}

To illustrate what the conditions on the parameters imply, we plot the phase diagram in the normalized PD (i.e. $t=1$ and $s=0$) in Figure \ref{fig:regionplotPDd65}. As expected, we see that increasing $r$ (the parameter controlling how much the players earn when in $(C, C)$) increases the number of possible solutions. Decreasing $p$ has a similar effect, but there is a critical value for $r$ that dictates whether it is possible to reach a region in which the GT or WSLS strategy pairs exist by decreasing $p$ or not. We show how the regions change as a function of $\delta$ in Appendix \ref{app:phasediagramdelta}.

The critical conditions in Table \ref{tab:absorbingPD} indicate when changes in the BRN  lead to the appearance or disappearance of equilibria. In Table \ref{tab:criticalconPD} we show all critical conditions at which changes in the BRN occur. These results give two extremes in the resolution of the phase diagram: the first is the lowest resolution and the second is the highest. There is an intermediate resolution in which we identify the critical conditions at which changes in the BRN graph lead to changes in the basin of attraction of the equilibria. We leave the study of the basins of attraction of the absorbing states for future work.  

\begin{table}[!ht]
\centering
\vspace{0.1cm}

\begin{tabular}{l c} \toprule
\textbf{Edge}   &\textbf{Conditions}  \\ \midrule

$0\rightarrow 0$ &  None\\ 
$1\rightarrow 0$ &  $\delta + r < 1 + \delta p$\\
$1\rightarrow 1$ &  $\delta + r > 1 + \delta p$\\
$2\rightarrow 0$ &  None\\
$3\rightarrow 0$ &  None\\
$4\rightarrow 0$ &  $\delta < p(1+\delta)$\\
$4\rightarrow 13$ &  $\delta > p(1 + \delta)$\\
$5\rightarrow 0$ &  $(p+r < 1 \wedge \delta < p(1+\delta))\vee (p+r \geq 1 \wedge \delta + r < p(1+\delta))$\\
$5\rightarrow 3$ &  $2r>1 \wedge p+r>1\wedge \frac{1-r}{1-p}<\delta<\frac{p}{r}$\\
$5\rightarrow 12$ &  $p+r < 1 \wedge \delta > p(1+\delta)\wedge r(1+\delta) < 1$\\
$5\rightarrow 15$ &  $(p+r \leq 1 \wedge r(1+\delta) \geq 1)\vee (p+r > 1 \wedge \delta r > p)$\\
$6\rightarrow 0$ &  $\delta < p$\\$7\rightarrow 10$ &  $\delta > p$\\
$7\rightarrow 0$ &  $\delta < p$\\$8\rightarrow 9$ &  $\delta > p$\\
$8\rightarrow 0$ &  None\\
$9\rightarrow 1$ &  $(1+\delta)r < 1 + \delta p$\\
$9\rightarrow 9$ &  $(1+\delta)r > 1 + \delta p$\\
$10\rightarrow 0$ &  None\\
$11\rightarrow 0$ &  None\\
$12\rightarrow 0$ &  None\\
$13\rightarrow 0$ &  $(1+\delta)r < 1 + \delta p$\\
$13\rightarrow 11$ &  $(1+\delta)r > 1 + \delta p$\\
$14\rightarrow 0$ &  None\\
$15\rightarrow 0$ &  None\\
 \bottomrule
\end{tabular}
\caption{Critical conditions for the existence of edges in the pure best-response graph of the prisoner's dilemma.}\label{tab:criticalconPD}
\end{table} 

\begin{figure}
\centering
\includegraphics[scale=0.5]{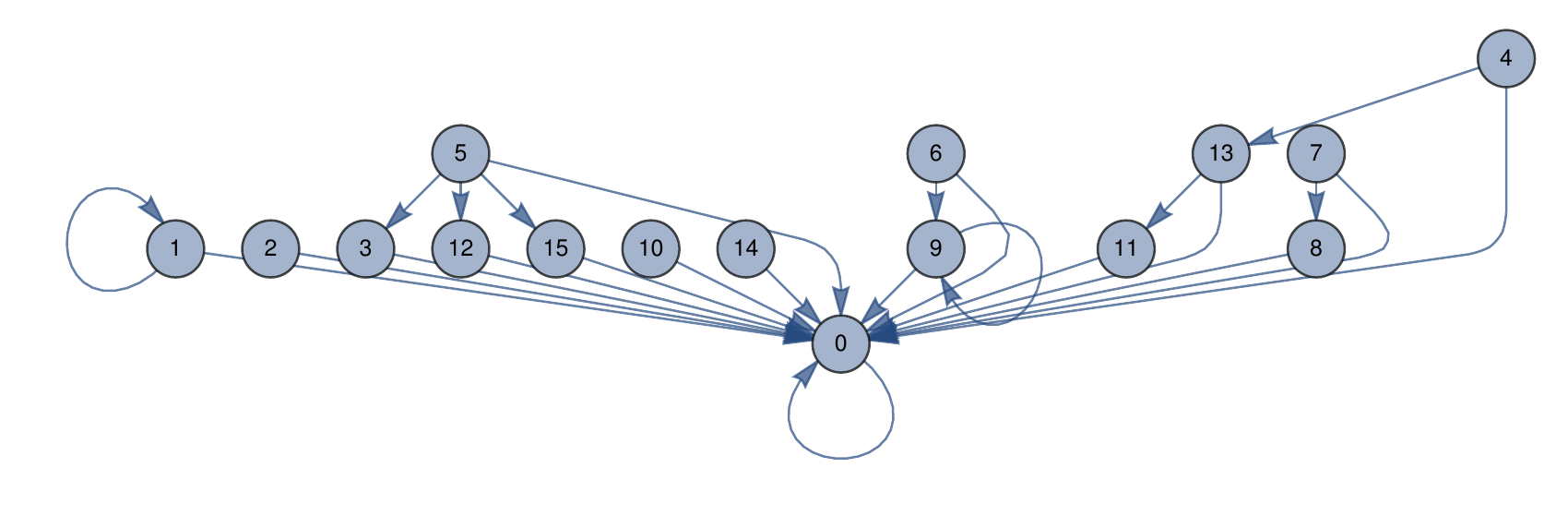}
\caption{All possible edges in the BRN of the prisoner's dilemma. See Table \ref{tab:criticalconPD} for conditions for the existence of edges.}
\label{fig:pBRgraphPD}
\end{figure}

From Table \ref{tab:criticalconPD} we can infer the number of BRN that are possible. Since some conditions are repeated (e.g. the conditions for nodes 9 and 13) we find that there are 12 distinct BRN (i.e. there are twelve regions in the phase space). We show the phase diagram for these 12 graphs in the left panel of Figure \ref{fig:regionplotPDd65}. In Figure \ref{fig:pBRgraphPD} we show an example of such a graph in which we see all three of the possible equilibria as self-loops for nodes 0, 1 and 9. 

\subsection{Stag Hunt}
\label{sec:resultsSH}
In the case of the stag hunt game, the number of equilibria increases significantly, which means that a simple assessment of the resulting phase diagram is not feasible. We see in Table \ref{tab:absorbingSH} that the SH gives rise to a total of 16 possible absorbing states of which half are symmetric and half are non-symmetric. The non-symmetric absorbing states are all of the the type CM defined in Section \ref{sec:meth:act}. We also find that almost 70\% of the absorbing states are conducive to cooperation. 

\begin{figure}
\centering
\includegraphics[scale=1.8]{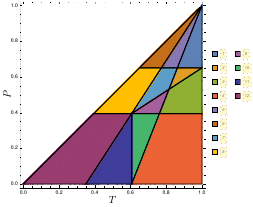}
\includegraphics[scale=1.65]{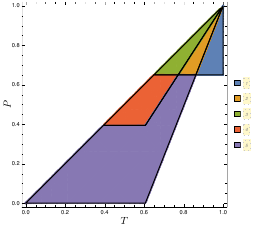}
\caption{Left: Phase diagram for the possible BRNs in the stag hunt Right: Phase diagram for the possible equilibria in the stag hunt. For both plots we have $s=0, r=1$ and $\delta=0.65$.}
\label{fig:regionplotSHd65}
\end{figure}

Table \ref{tab:criticalconSH} shows some similarities to Table \ref{tab:criticalconPD} with many critical conditions being identical when exchanging $t$ and $r$. This similarity is confirmed by the left panel of Figure \ref{fig:regionplotSHd65} where we see the same 12 regions appearing as in the left panel of Figure \ref{fig:regionplotPDd65}. There is, however, a difference in the regions that are relevant for changes in the set of absorbing states as seen in the right panel of Figure \ref{fig:regionplotSHd65}.  
 
\begin{figure}
\centering
\includegraphics[scale=0.5]{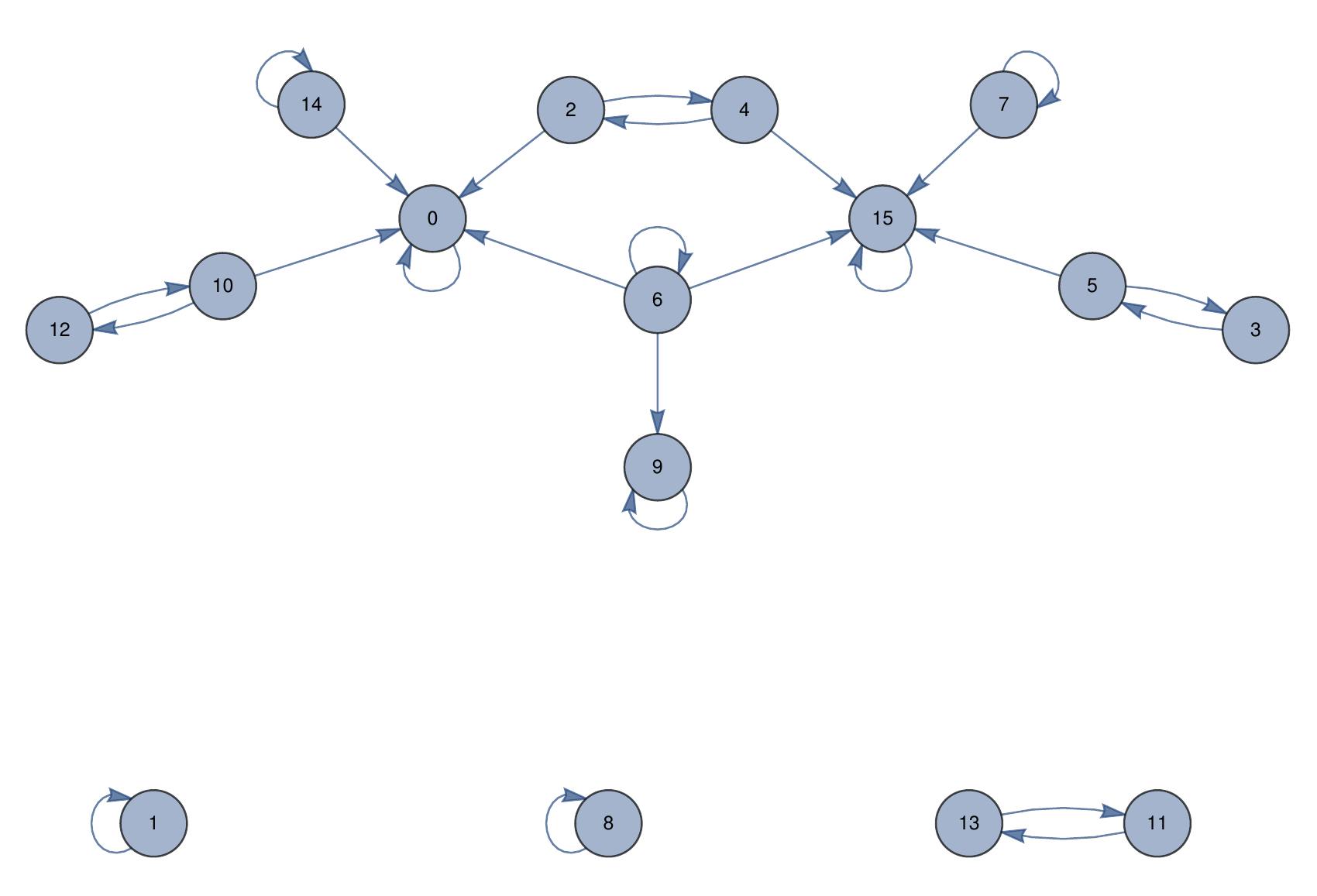}
\caption{All possible edges in the BRN of the stag hunt. See Table \ref{tab:criticalconSH} for conditions for the existence of edges.}
\label{fig:pBRgraphSH}
\end{figure}

In Figure \ref{fig:pBRgraphSH} we see that the network of possible edges in the BRN also differs from the PD case. More specifically, we see that the maximum in degree in the SH is 5 while being 16 in the PD. This suggest that the basins of attraction of the absorbing states are more evenly distributed in the SH than in the PD (the IBRN shown in Appendix \ref{app:IBRN} support this hypothesis).  

\begin{table}[!ht]
    \centering
    \begin{tabular}{c|c|c|c}\toprule
        \textbf{Policy-pair} & \textbf{Type} & \textbf{Conditions} & \textbf{Region}\\
        \hline
        $0$ & S, NCC & Always & $(0)$\\
        $17$ & S, NCC & Always & $(0)$\\
        $102$ & S, NCC & $p+\delta p > s + \delta r \land r + \delta p > t + \delta r$ & $(36)$\\
        $119$ & S, CC & $p+\delta s > s + \delta r$ & $(53)$\\
        $136$ & S, CC & Always & $(0)$\\
        $153$ & S, CC & Always & $(0)$\\
        $238$ & S, CC & $r + \delta p > t + \delta t$ & $(172)$\\
        $255$ & S, CC & Always & $(0)$\\
        \hline
        $36$ & CM, NCC & $p+\delta p > s + \delta r \land r + \delta p > t + \delta r$ & $(36)$\\
        $53$ & CM, CC & $p+\delta s > s + \delta r$ & $(53)$\\
        $66$ & CM, NCC  & $p+\delta p > s + \delta r \land r + \delta p > t + \delta r$ & $(36)$\\
        $83$ & CM, CC & $p+\delta s > s + \delta r$ & $(53)$\\
        $172$ & CM, CC & $r + \delta p > t + \delta t$ & $(172)$\\
        $189$ & CM, CC & Always & $(0)$\\
        $202$ & CM, CC & $r + \delta p > t + \delta t$ & $(172)$\\
        $219$ & CM, CC & Always & $(0)$\\
         \bottomrule
    \end{tabular}
    \caption{All strategy pairs solving the Bellman equations self-consistently for the stag hunt.}
    \label{tab:absorbingSH}
\end{table}

\begin{table}[!ht]
\centering
\vspace{0.1cm}

\begin{tabular}{l c} \toprule
\textbf{Edge}   &\textbf{Conditions}  \\ \midrule

$0\rightarrow 0$ &  None\\ 
$1\rightarrow 1$ &  None\\
$2\rightarrow 0$ &  $\delta + t > 1+ \delta p$\\
$2\rightarrow 4$ &  $\delta + t < 1+ \delta p$\\
$3\rightarrow 5$ &  None\\
$4\rightarrow 2$ &  $\delta < \delta(1+p)$\\
$4\rightarrow 15$ &  $\delta < \delta(1+p)$\\
$5\rightarrow 3$ &  $\delta<p$\\
$5\rightarrow 15$ &  $\delta>p$\\
$6\rightarrow 0$ &  $2t>1 \wedge p+t>1\wedge \frac{1-t}{1-p}<\delta<\frac{p}{t}$\\
$6\rightarrow 6$ &  $(p+t < 1 \wedge \delta < p(1+\delta))\vee p+t\geq 1 \wedge \delta+t<1+\delta p$\\
$6\rightarrow 9$ &  $(p+t \leq 1 \wedge 1 < t(1+\delta))\vee p+t> 1 \wedge \delta t< p$\\
$6\rightarrow 15$ &  $p+t<1 \wedge \delta>p(1+\delta)\wedge t(1+\delta)<1$\\
$7\rightarrow 7$ &  $\delta < p$\\
$7\rightarrow 15$ &  $\delta > p$\\
$8\rightarrow 8$ &  None\\
$9\rightarrow 9$ &  None\\
$10\rightarrow 0$ &  $(1+\delta)t > 1 + \delta p$\\
$10\rightarrow 12$ &  $(1+\delta)t < 1 + \delta p$\\
$11\rightarrow 13$ &  None\\
$12\rightarrow 10$ &  None\\
$13\rightarrow 11$ &  $(1+\delta)r < 1 + \delta p$\\
$14\rightarrow 0$ &  $(1+\delta)t > 1 + \delta p$\\
$14\rightarrow 14$ &  $(1+\delta)t < 1 + \delta p$\\
$15\rightarrow 0$ &  None\\
 \bottomrule
\end{tabular}
\caption{Critical conditions for the existence of edges in the BRN of the stag hunt.}\label{tab:criticalconSH}
\end{table} 

\subsection{Hawk-Dove}
\label{sec:resultsHD}

For the hawk-dove game, the number of equilibria increases significantly, is similar to that of the SH. We see in Table \ref{tab:absorbingHD} that the HD gives rise to a total of 17 possible absorbing states of which 16 are not conducive to cooperation. The absorbing state that is CC is the symmetric WSLS, which is also an absorbing state for both of the other games. Half of the remaining absorbing states are C and the other half are CS. 

\begin{figure}
\centering
\includegraphics[scale=1.8]{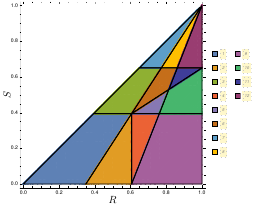}
\caption{Phase diagram for the possible networks in the hawk-dove game with $s=0, r=1$ and $\delta=0.65$.}
\label{fig:regionplotHDd65}
\end{figure}

Table \ref{tab:criticalconHD} again shows some similarities to Table \ref{tab:criticalconPD} with many critical conditions being identical when exchanging $s$ and $p$. Once again we find the 12 regions seen in the left panel of Figure \ref{fig:regionplotPDd65} in the phase diagram in Figure \ref{fig:regionplotHDd65}. In this case, all regions are relevant for changes in the set of absorbing states.  
 
\begin{figure}
\centering
\includegraphics[scale=0.5]{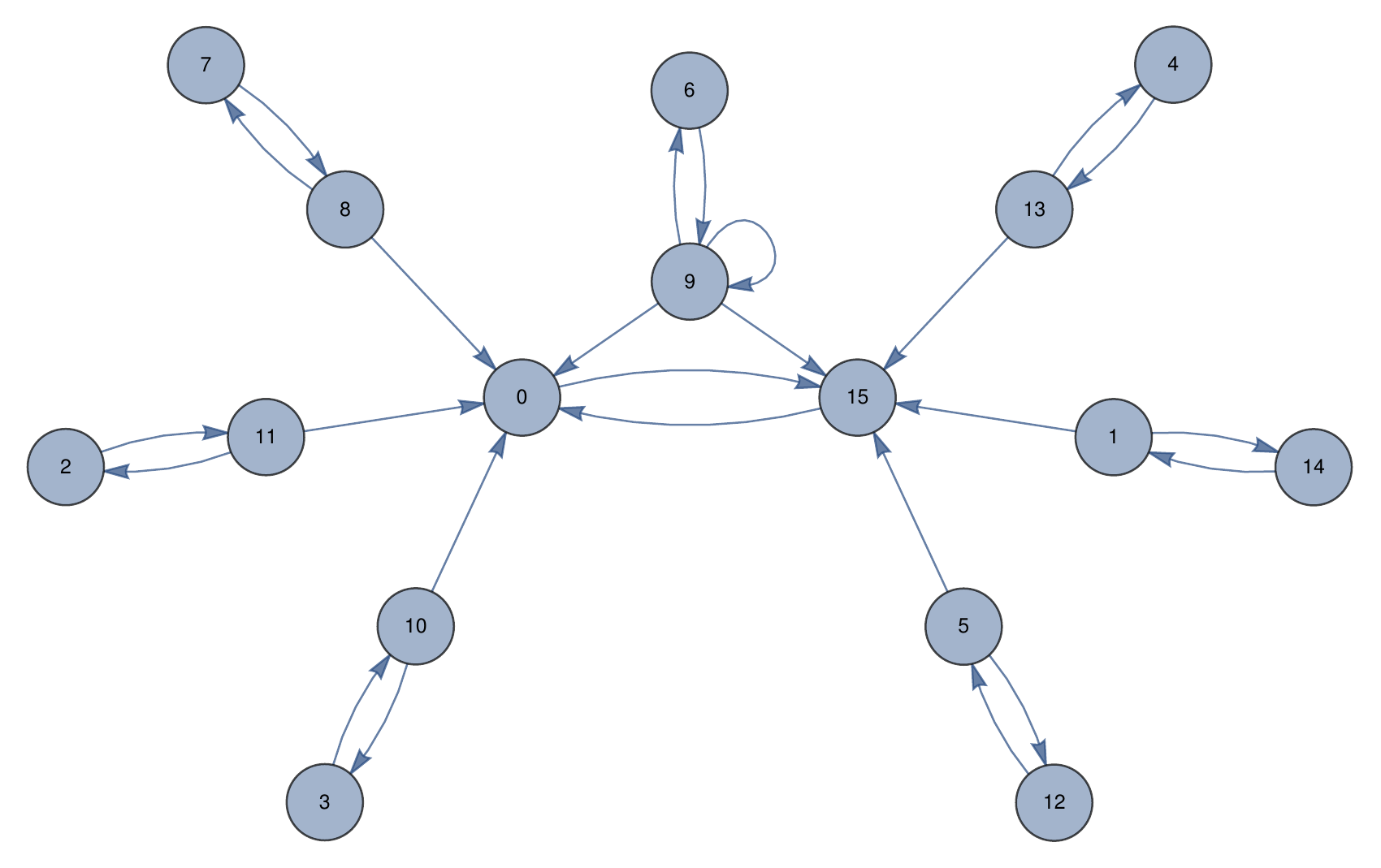}
\caption{All possible edges in the BRNs of the hawk dove game. See Table \ref{tab:criticalconHD} for conditions for the existence of edges.}
\label{fig:pBRgraphHD}
\end{figure}

In Figure \ref{fig:pBRgraphHD} we see that the network of possible edges in the BRN for the HD game is similar to the network for the SH. More specifically, we see that the maximum in degree in the HD is also 5. Interestingly, the nodes realizing this maximum are the same for both games (nodes 0 and 15, i.e. the All-D and all cooperate strategies) with the difference that they have a self-loop in the SH, but do not have a self-loop in the HD game. 

\begin{table}[!ht]
    \centering
    \begin{tabular}{c|c|c|c}\toprule
        \textbf{Policy-pair} & \textbf{Type} & \textbf{Conditions} & \textbf{Region}\\
        \hline
        $153$ & S, CC & $p + \delta r > s + \delta p \land r + \delta r > t + \delta p$ & $(153)$\\
        \hline
        $15$ & C, NCC & Always & $(15)$\\
        $30$ & C, NCC & $t + \delta s > r + \delta t$ & $(30)$\\
        $105$ & C, NCC & $s + \delta s > p + \delta t \land t + \delta s > r + \delta t$ & $(105)$\\
        $120$ & C, NCC & $s + \delta s > p + \delta t$ & $(120)$\\
        $135$ & C, NCC & $s + \delta s > p + \delta t$ & $(120)$\\
        $150$ & C, NCC & $s + \delta s > p + \delta t \land t + \delta s > r + \delta t$ & $(105)$\\
        $225$ & C, NCC & $t + \delta s > r + \delta t$ & $(30)$\\
        $240$ & C, NCC & Always  & $(15)$\\
        \hline
        $43$ & CS, NCC & $s + \delta t > p + \delta t$ & $(43)$\\
        $58$ & CS, NCC & $s + \delta t > p + \delta t$ & $(43)$\\
        $77$ & CS, NCC & $t + \delta s > r + \delta r$ & $(77)$\\
        $92$ & CS, NCC & $t + \delta s > r + \delta r$ & $(77)$\\
        $163$ & CS, NCC & $s + \delta t > p + \delta t$ & $(43)$\\
        $178$ & CS, NCC & $s + \delta t > p + \delta t$ & $(43)$\\
        $197$ & CS, NCC & $t + \delta s > r + \delta r$ & $(77)$\\
        $212$ & CS, NCC & $t + \delta s > r + \delta r$ & $(77)$\\
        \bottomrule
    \end{tabular}
    \caption{All strategy pairs solving the Bellman equations self-consistently for the hawk-dove game.}
    \label{tab:absorbingHD}
\end{table}

\begin{table}[!ht]
\centering
\vspace{0.1cm}

\begin{tabular}{l c} \toprule
\textbf{Edge}   &\textbf{Conditions}  \\ \midrule

$0\rightarrow 15$ &  None\\ 
$1\rightarrow 14$ &  $\delta + r < 1+\delta s$\\
$1\rightarrow 15$ &  $\delta + r > 1+\delta s$\\
$2\rightarrow 11$ &  None\\
$3\rightarrow 10$ &  None\\
$4\rightarrow 13$ &  None\\
$5\rightarrow 12$ &  $r(1+\delta) < 1+ \delta s$\\
$5\rightarrow 15$ &  $r(1+\delta) > 1+ \delta s$\\
$6\rightarrow 9$ &  None\\
$7\rightarrow 8$ &  None\\
$8\rightarrow 0$ &  $\delta > s (1+\delta)$\\
$8\rightarrow 7$ &  $\delta < s (1+\delta)$\\
$9\rightarrow 0$ &  $r+s < 1 \wedge \delta < s(1+\delta) \wedge r(1+\delta)< 1$\\
$9\rightarrow 6$ &  $(r+s \leq 1 \wedge \delta < s(1+\delta)) \vee (r+s >1 \wedge \delta+r < 1 +\delta s)$\\
$9\rightarrow 9$ &  $(r+s < 1 \wedge 1 > r(1+\delta)) \vee (r+s \geq 1 \wedge \delta r < s)$\\
$9\rightarrow 15$ &  $\frac{1-r}{1-s}<\delta<\frac{s}{r}\wedge (2s>1 \vee r+s>1)$\\
$10\rightarrow 0$ &  $\delta>s$\\
$10\rightarrow 3$ &  $\delta<s$\\
$11\rightarrow 0$ &  $\delta>s$\\
$11\rightarrow 2$ &  $\delta<s$\\
$12\rightarrow 5$ &  None\\
$13\rightarrow 4$ &  $(1+\delta)r < 1 + \delta s$\\
$13\rightarrow 15$ &  $(1+\delta)t > 1 + \delta s$\\
$14\rightarrow 1$ &  None\\
$15\rightarrow 0$ &  None\\
 \bottomrule
\end{tabular}
\caption{Critical conditions for the existence of edges in the BRN of the hawk dove game.}\label{tab:criticalconHD}
\end{table} 

\section{Sample Batch Q-learning}
\label{sec:numerics}
The learning dynamics we have discussed thus far have been framed predominantly in terms of mutual best-response dynamics. In this section, we show that these dynamics are the limiting dynamics of sample batch Q-learning when taking the infinite batch size limit. We will use a simple sample batch algorithm so that it is simple to see that the infinite batch size limit gives rise to the mutual best-response dynamics, but more sophisticated sample batch algorithms exists which make better use of the information gathered during a batch (see \cite{Barfuss2019, Barfuss2020}).

The sample batch algorithm splits the time horizon into batches of size $K$ and keeps track of two sets of Q-values. The first set, denoted by $Q^{\text{act}}$, will be used to determine the actions that the algorithm takes during a batch. These Q-values are kept constant during the batch. The second set, denoted by $Q^{\text{val}}$, will be used to learn during the batch. At the end of the batch the first set is replaced by the second and the process is repeated. 

In order to ensure that all states are visited infinitely often in the infinite batch size limit, the algorithm uses an $\epsilon$-greedy action-selection mechanism with decreasing exploration rate during each batch. This means that it selects the action with the maximum Q-value with probability $1-\epsilon$ and selects an action uniformly at random with probability $\epsilon$. For each batch we enumerate the periods by $k\in \{1, 2, \ldots, K\}$ and use an exploration rate of 
\begin{equation}
\epsilon_{k} = \frac{1}{\sqrt{k}}
\end{equation}
in the $k^{\text{th}}$ period. This ensure that the states that are only visited when both agents explore (which occurs with probability $\sim\epsilon^{2}$) are also visited infinitely often in the infinite batch size limit. 

During the batch the second set of Q-values is updated using
\begin{equation}
Q^{\text{val}}_{a, \boldsymbol{\sigma}, \sigma_{a}}(k+1) = (1-\tilde{\alpha}(k))Q^{\text{val}}_{a, \boldsymbol{\sigma}, \sigma_{a}}(k) + \tilde{\alpha}(k)\big[r_{a}(\boldsymbol{\sigma}(k+1)) + \delta \max_{i\in C, D} \{Q^{\text{val}}_{a, \boldsymbol{\sigma}}(k+1), i\}\big]
\end{equation}
with 
\begin{equation}
\tilde{\alpha}(k)=
\begin{cases}
\frac{1}{k_{\boldsymbol{\sigma}, \sigma_{a}}(k)+1}\quad &\text{if} \; (\boldsymbol{\sigma}, \sigma_{a}) = (\boldsymbol{\sigma}(k), \sigma_{a}(k))\\
0  &\text{otherwise}
\end{cases},
\end{equation}
where $k_{\boldsymbol{\sigma}, \sigma_{a}}(k)$ is a local time, given by the number of times that the state-action pair $(\boldsymbol{\sigma}, \sigma_{a})$ has been visited up to period $k$. At the end of a batch the Q-values used for determining the actions are updated as
\begin{equation}
Q^{\text{act}}_{a, \boldsymbol{\sigma}, \sigma_{a}} = Q^{\text{val}}_{a, \boldsymbol{\sigma}, \sigma_{a}}(K), 
\end{equation}

The Q-value update during the batch defined as above satisfies the conditions for convergence given in \cite{Watkins1992}, since the environment in which the updates are made is stationary. In the infinite batch size limit the Q-values are updated to the solution of the Bellman equations given the opponents strategy as calculated in Section \ref{sec:selfcons}. 

\section{Conclusion}
\label{sec:conclusion}
We have developed a computational method to identify absorbing states for two-player, two-action repeated games played by Q-learners with memory one. For the PD we have shown that the three symmetric solutions to the Bellman equations identified by \cite{Usui2021} (All-D, GT and WSLS) are the only pure strategy pair solutions possible. We also find the absorbing states and the conditions for their existence for the stag hunt and hawk-dove games. 

Furthermore, we have developed a graphical representation of the mutual best-response dynamics and have shown that these are the limiting dynamics of sample batch Q-learning  with infinite batch size. The representation shows that limit cycles become possible. This is an artifact of the Q-value updates taking place simultaneously in Sample Batch Q-learning, as opposed to the sequential Q-value updates in the algorithm proposed in \cite{Usui2021}.

The absorbing states identified using these methods can be used to identify the learning goals for Q-learning algorithms with memory in these settings. The advantage of using these as learning goals is that the mutual best-response dynamics could provide a method for establishing convergence results for sample batch Q-learning algorithms with finite but increasing batch size. The most immediate application of the method is to study the basins of attractions of the absorbing states (this would have to include the basins of attraction of all limit cycles). By classifying the limiting dynamics (absorbing states and limit cycles) into either CC or NCC we can identify the likelihood of cooperating for each region of the parameter space and study how the likelihood of cooperation changes as a function of the model parameters. 

The methods developed here lend themselves to extensions such as considering two period memory algorithms, considering three action games or placing alternative conditions on the model parameters. Another direction is the extension of the results to off-policy algorithms, or other temporal difference learning algorithms. The method can also be extended to study how the networks change when considering a smaller learning rate, i.e., where the algorithms move in the direction of the best-response instead of playing the best-response immediately. Another extension would be to study the effect of noise due to using large but finite batch sizes. The previous two extensions are being investigated in parallel work \cite{Barfuss2022}. Finally, the method may be extended to the sequential move setting by constructing a directed network with two types of edges: the first dictating how the strategy pair changes when player one acts and the second dictating how the strategy pair changes when player two acts. The sequential best-response dynamics in this network will follow paths of alternating edge types.  

Our results, together with the extensions outlined above provide an avenue for understanding the dynamics of multi-agent reinforcement learning algorithms in simple settings. The insights gained through this understanding may also be applicable to more complex environments and may be used in designing algorithms that learn successfully in such settings.

\appendix

\section{Phase diagram as a function of the discount factor}     
\label{app:phasediagramdelta}   
\begin{figure}
\centering
\includegraphics[scale=1.7]{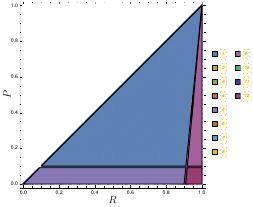}
\includegraphics[scale=1.7]{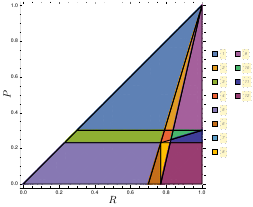}
\includegraphics[scale=1.7]{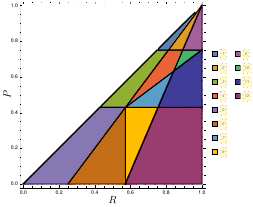}
\includegraphics[scale=1.7]{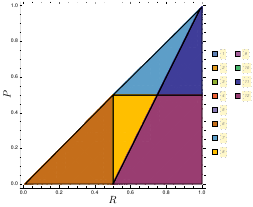}
\caption{Phase diagram for the possible BRNs in the prisoner's dilemma for different values of $\delta$: 0.1 (Top Left), 0.3 (Top Right), 0.75 (Bottom Left) and 0.99 (Bottom Right).}
\label{fig:regionplotPDasDelta}
\end{figure}

In Figure \ref{fig:regionplotPDasDelta} we plot the 12 regions of the phase diagram for the PD for different values of $\delta$. This shows how the regions change as we vary the discount factor. The plots for the SH and HD games will show the same behavior with the difference that the regions represent different graphs. In particular, the regions associated with the existence of specific absorbing states are different. 

The plots show that an increase in $\delta$ leads to an increase in fraction of the phase space occupied by the region where all three possible equilibria exist in the PD. For the SH and HD in contrast, we find that regions in which certain combinations of absorbing states exist disappear as we increase $\delta$.  

\section{Examples of IBRN}
\label{app:IBRN}
\begin{figure}
\centering
\includegraphics[scale=0.35]{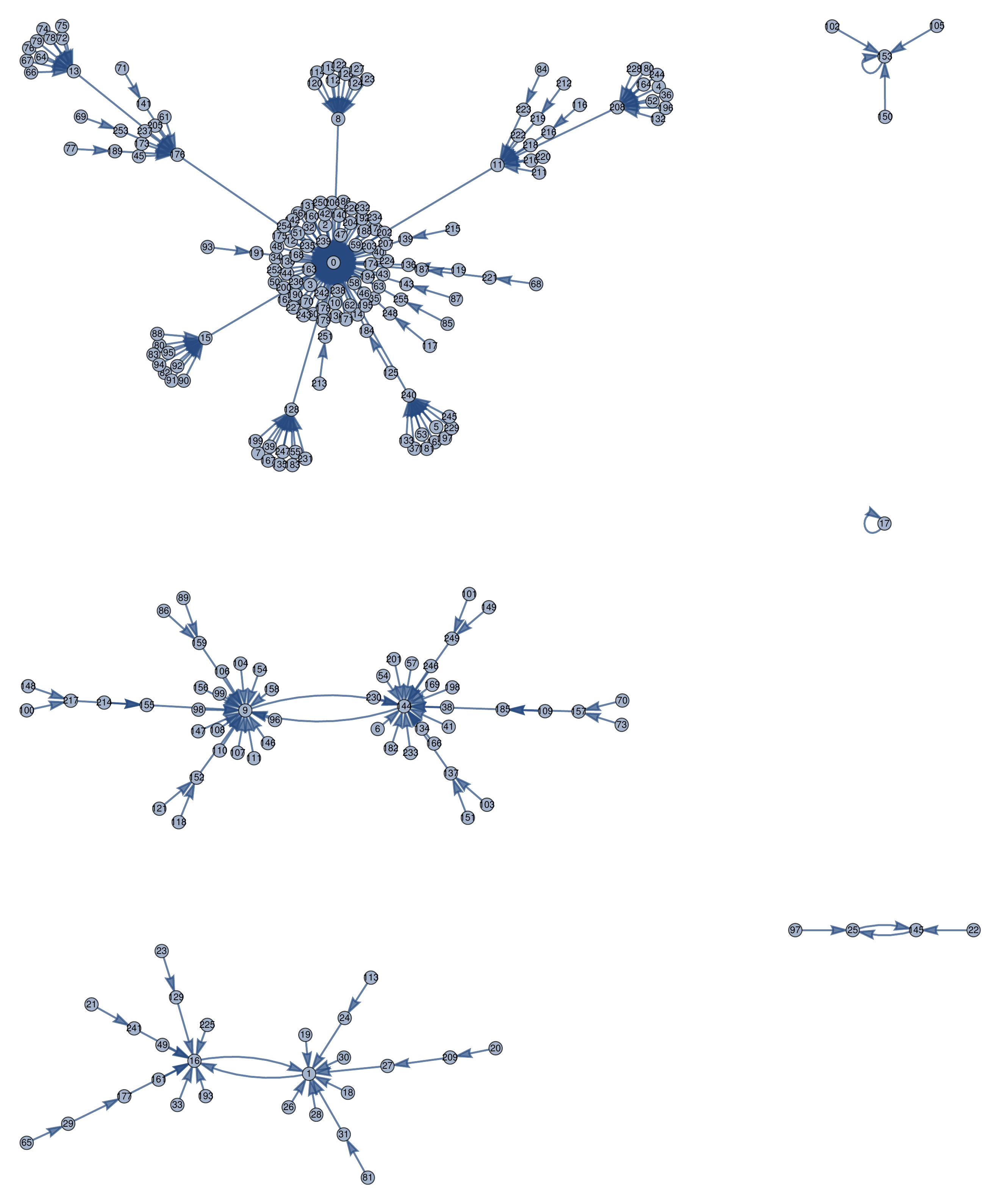}
\caption{IBRN in the region 12 of the PD.}
\label{fig:IBRNPD12}
\end{figure}

\begin{figure}
\centering
\includegraphics[scale=0.35]{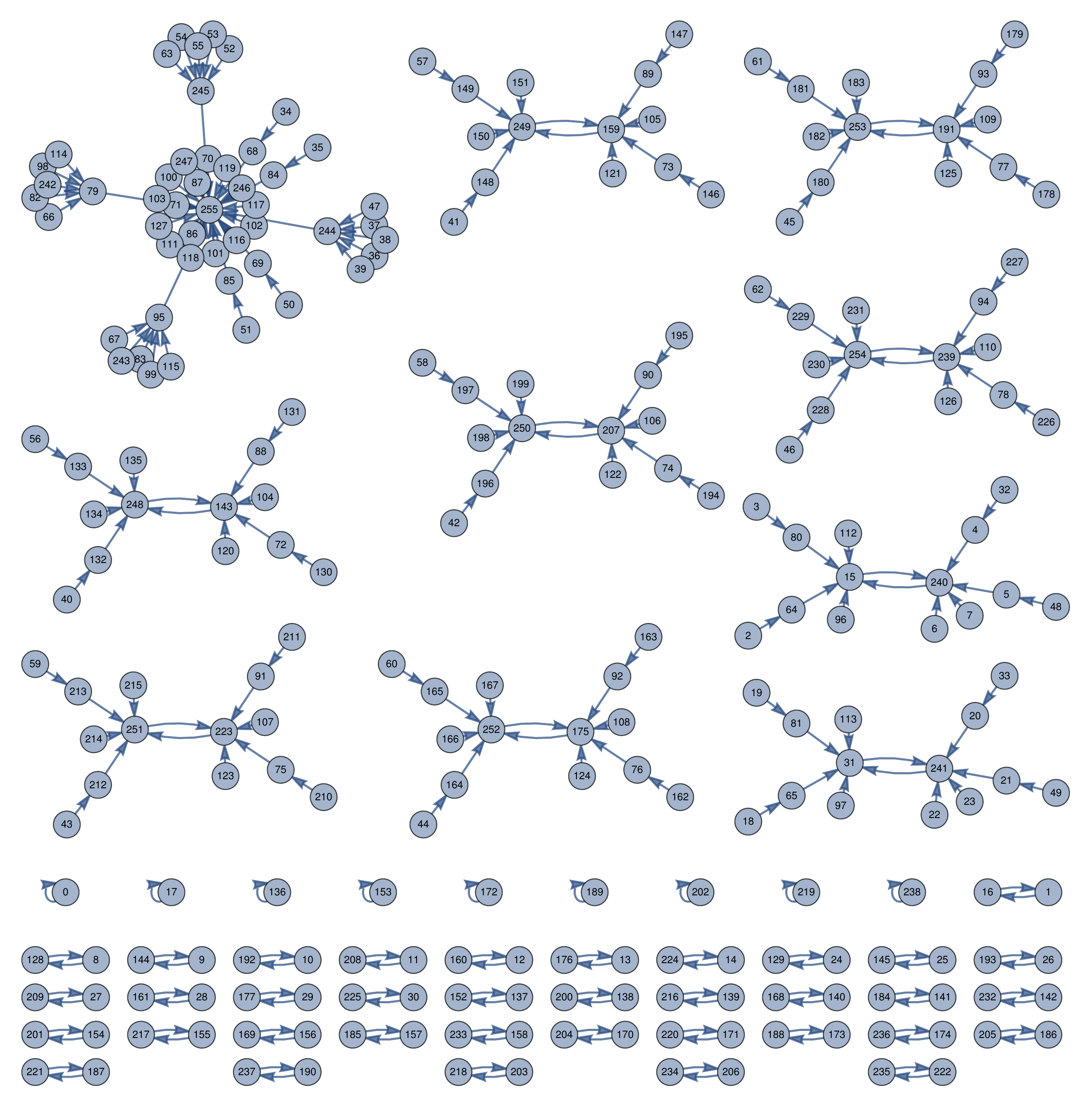}
\caption{IBRN in the region 12 of the SH.}
\label{fig:IBRNSH12}
\end{figure}

\begin{figure}
\centering
\includegraphics[scale=0.35]{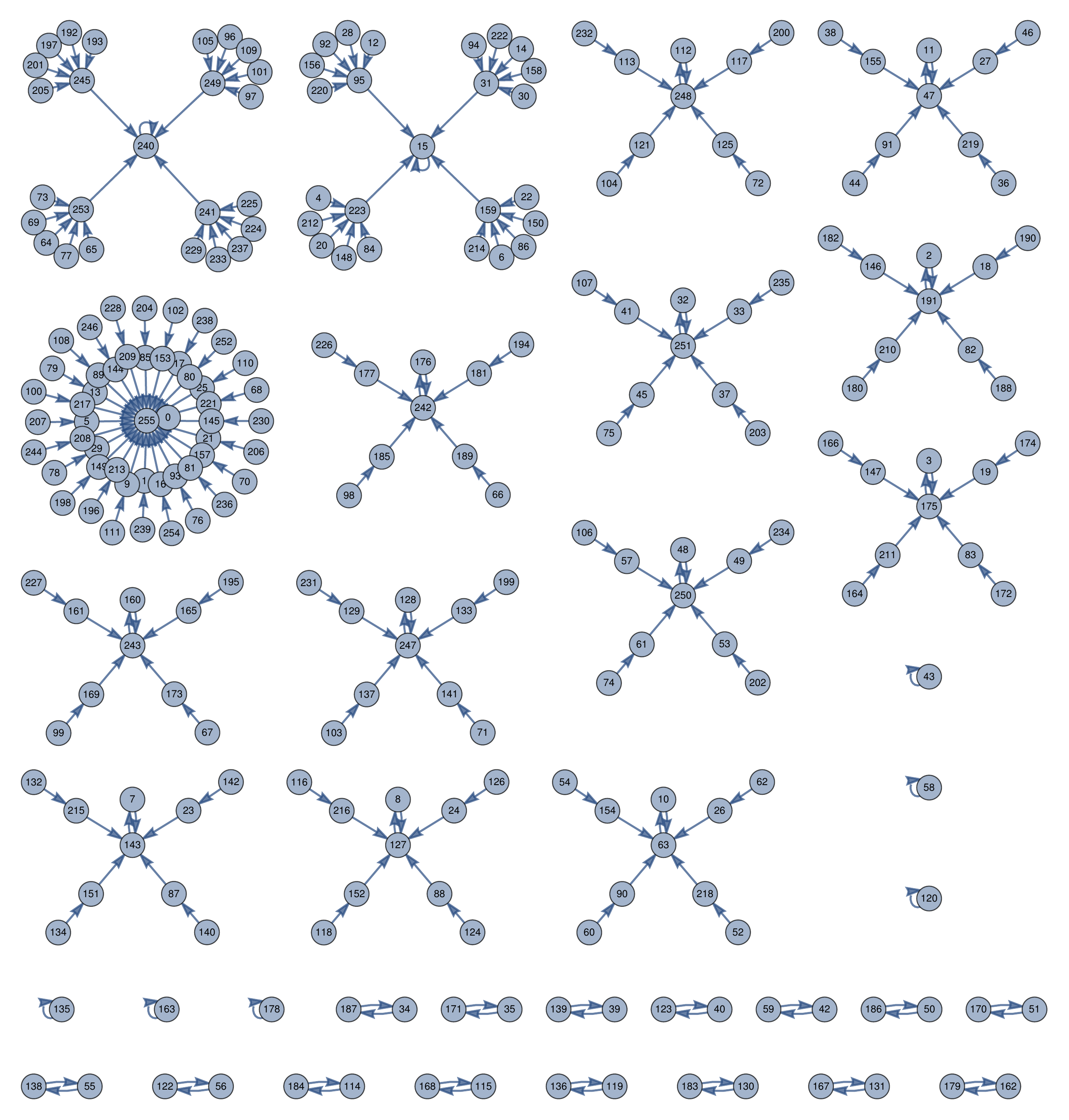}
\caption{IBRN in the region 12 of the HD game.}
\label{fig:IBRNHD12}
\end{figure}

In Figure \ref{fig:IBRNPD12}, Figure \ref{fig:IBRNSH12} and Figure \ref{fig:IBRNHD12}, we show the IBRN in the $12^{\text{th}}$ region of the PD, SH and HD respectively. We indeed see that the sizes of the basins of attraction for the absorbing states are more evenly distributed in the SH and HD than in the PD.

\appendix


\newpage
\begin{thebibliography}{99}

\bibitem{Axelrod1981}
Axelrod, R., and Hamilton, W.D. (1981) \emph{The evolution of cooperation}, Science 211 (4489) 1390--1396 .

\bibitem{Barfuss2020}
Barfuss, W.  (2020) \emph{Reinforcement learning dynamics in the infinite memory limit} in Proceedings of the 19th International Conference on Autonomous Agents and MultiAgent Systems (International Foundation for Autonomous Agents and Multiagent Systems, 2020)

\bibitem{Barfuss2019}
Barfuss, W., Donges, J.F., and Kurths, J. (2019) \emph{Deterministic limit of temporal difference reinforcement learning for stochastic games}, Phys. Rev. E 99, 043305. 

\bibitem{Barfuss2022}
Barfuss, W. and Meylahn, J.M. (2022) \emph{Intrinsic fluctuations of reinforcement learning promote cooperation}, to be submitted. 

\bibitem{Busoniu2008}
Busoniu, L., Babu\v{s}ka, R., and De Schutter, B. (2008). \emph{A comprehensive survey of multiagent reinforcement learning}. IEEE Transactions on Systems, Man, and Cybernetics, Part C (Applications and Reviews), 38(2), 156--172.

\bibitem{Busoniu2010}
Bu\c{s}oniu, L., Babu\v{s}ka, R., and De Schutter, B. (2010). \emph{Multi-agent reinforcement learning: An overview}. Innovations in multi-agent systems and applications-1, 183--221.

\bibitem{Calvano2020}
Calvano, E., Calzolari, G., Denicol\'{o}, V., Harrington, J. and Pastorello, S. (2020) \emph{Protecting Consumers from high prices due to AI}, Science, 370(6520):1040--1042.

\bibitem{Calvano2020b}
Calvano, E., Calzolari, G., Denicol\'{o}, V. and Pastorello, S. (2020b) \emph{Artificial Intelligence, Algorithmic Pricing and Collusion}, American Economic Review 110(10), 3267-97.

\bibitem{Fudenberg1991}
Fudenberg, Drew, and Tirole, J. (1991). \emph{Game theory}. MIT press.

\bibitem{Hennes2008}
Hennes, D., Tuyls, K., and Rauterberg, M. (2008). \emph{Formalizing multi-state learning dynamics}. In 2008 IEEE/WIC/ACM International Conference on Web Intelligence and Intelligent Agent Technology 2,266--272.

\bibitem{Hernandez2019}
Hernandez-Leal, P., Kartal, B., and Taylor, M. E. (2019). \emph{A survey and critique of multiagent deep reinforcement learning}. Autonomous Agents and Multi-Agent Systems, 33(6), 750--797.

\bibitem{Jaakkola1994}
Jaakkola, T., Jordan, M. I., and Singh, S. P. (1994). \emph{On the convergence of stochastic iterative dynamic programming algorithms}. Neural computation, 6(6), 1185--1201.

\bibitem{Klein2021}
Klein, T. (2021). \emph{Autonomous Algorithmic Collusion: Q-Learning Under Sequential Pricing}, The RAND Journal of Economics 52(3), 538-558.

\bibitem{Macy2002}
Macy, M. W., and Flache, A. (2002). \emph{Learning dynamics in social dilemmas}. Proceedings of the National Academy of Sciences, 99(suppl 3), 7229--7236.

\bibitem{Meylahn2022}
Meylahn, J. M., and V. den Boer, A. (2022). \emph{Learning to collude in a pricing duopoly}, Manufacturing \& Service Operations Management, 0(0).

\bibitem{Palmer2020}
Palmer, G. \emph{Independent learning approaches: Overcoming multi-agent learning pathologies in team-games}, PhD Thesis (2020).

\bibitem{Panait2005}
Panait, L., and Luke, S. (2005). \emph{Cooperative multi-agent learning: The state of the art}. Autonomous agents and multi-agent systems, 11(3), 387--434.

\bibitem{Rabin1993}
Rabin, M. (1993). \emph{Incorporating fairness into game theory and economics}. The American economic review, 1281--1302.

\bibitem{Rapoport1965}
Rapoport, A., Chammah, A. M., and Orwant, C. J. (1965). \emph{Prisoner's dilemma: A study in conflict and cooperation} (Vol. 165). University of Michigan press.

\bibitem{Rousseau1754}
Rousseau, J-J. (1754) \emph{Discours sur l'origine et les fondements de l'in\'{e}galit\'{e} parmi les hommes}, Paris : Bordas.

\bibitem{Smith1973}
Smith, JMPGR and Price, G. R. (1973) \emph{The logic of animal conflict}. Nature, 5427(246), 15--18.

\bibitem{Usui2021}
Usui,  Y. and Ueda, M. (2021) \emph{Symmetric equilibrium of multi-agent reinforcement learning in repeated prisoner's dilemma}, Applied Mathematics and Computation 409, 126370.

\bibitem{Waltman2008}
Waltman, L., and Kaymak, U. (2008). \emph{Q-learning agents in a Cournot oligopoly model}. Journal of Economic Dynamics and Control, 32(10), 3275--3293.

\bibitem{Watkins1992}
Watkins, C. J., and Dayan, P. (1992). \emph{Q-learning}. Machine learning, 8(3-4), 279--292.

\bibitem{Wunder2010}
Wunder, M., Littman, M. L., and Babes, M. (2010, January). \emph{Classes of multiagent {Q}-learning dynamics with epsilon-greedy exploration}. In ICML.

\bibitem{Zhang2021}
Zhang, K., Yang, Z., and Ba\c{s}ar, T. (2021). \emph{Multi-agent reinforcement learning: A selective overview of theories and algorithms}. Handbook of Reinforcement Learning and Control, 321--384.

\end{thebibliography}
\end{document}